\documentclass[12pt,reqno,papersize]{amsart}
\usepackage{amsmath,amssymb}
\usepackage{graphicx,color,xcolor}
\usepackage{latexsym,amssymb}
\usepackage{mathrsfs}

\newtheorem{thm}{Theorem}[section]

\newtheorem{conjecture}[thm]{Conjecture}

\numberwithin{equation}{section}
\numberwithin{thm}{section}

\begin{document}
\begin{center}\large \bf 
Backward self-similar solutions for compressible 
\\ Navier-Stokes equations 
\end{center}

\vskip3mm

\footnote[0]
{


{\it Acknowledgements}: 
P. Germain was partially supported by the NSF grant DMS-1501019. 
T. Iwabuchi was supported by JSPS KAKENHI Grant 17H04824. 

{\it Addresses}: 
Pierre Germain, Courant Institute of Mathematical Sciences, New York University, 251 Mercer Street, New York, N.Y. 10012-1185, USA, {\tt pgermain@cims.nyu.edu}

Tsukasa Iwabuchi, Mathematical Institute, Tohoku University, 
Sendai City, 980-8578, JAPAN, 
{\tt t-iwabuchi@tohoku.ac.jp}

Tristan L\'eger, Courant Institute of Mathematical Sciences, New York University, 251 Mercer Street, New York, N.Y. 10012-1185, USA, {\tt tleger@cims.nyu.edu}
}


\begin{center}

Pierre GERMAIN, \quad Tsukasa IWABUCHI, \quad Tristan L\'EGER 

\vskip2mm


\end{center}

\vskip5mm

\begin{center}
\begin{minipage}{135mm}
\footnotesize
{\sc Abstract. } 
This article is devoted to backward self-similar blow up solutions of the compressible Navier-Stokes equations with radial symmetry. We show that such solutions cannot exist if they either satisfy an appropriate smallness condition, or have finite energy. Furthermore, numerical simulations lead us to the conjecture that such solutions do not exist. \end{minipage}
\end{center}



\section{Introduction}
\subsection{The model}
We study self-similar solutions of 
the following compressible Navier-Stokes equations 
in $\mathbb R^d$ with $d \geq 1$. 
\begin{equation}\label{eq:cNS_full}
\begin{cases}
\partial_t \rho + {\rm div \, } \big( \rho u \big) 
= 0 , 
& \quad t > 0, x \in \mathbb R^d, 
\\
\displaystyle 
\partial _t (\rho u) + 
{\rm div \,} \big( \rho u \otimes u \big) + \nabla \pi  
 = {\rm div } \, \tau,
& \quad t > 0, x \in \mathbb R^d, 
\\
\displaystyle 
\partial_t 
\Big[ \rho \Big( \frac{|u|^2}{2} + e \Big) \Big] 
+ {\rm div} 
\Big[ u \Big( \rho \Big(\frac{|u|^2}{2} + e \Big) 
      + \pi \Big) \Big]
- {\rm div } \, q = {\rm div \,} (\tau \cdot u) , 
& \quad t > 0, x \in \mathbb R^d, 
\end{cases}
\end{equation}
where $\rho(t,x)$ is the density of the fluid, $u(t,x)$ its velocity, $e(t,x)$ its internal energy, $\pi(t,x)$ its pressure, $\tau(t,x)$ its stress tensor, and finally $q(t,x)$ its internal energy flux. The fluid will furthermore be described by its temperature $\theta(t,x)$. 

We assume the following constitutive relations:
\begin{itemize}
\item Joule's first law: $$e = C_V \theta,$$ where $C_V > 0$ is the heat constant.
\item Ideal gas law: $$\pi = \rho R \theta,$$ where $R>0$ is the ideal gas constant.
\item Newtonian fluid: this implies 
$$\tau := \lambda {\rm div \,} u \, {\rm Id} 
 + 2 \mu D(u), 
\quad 
 D(u) = \frac{\nabla u + (\nabla u)^T}{2} , 
\quad \nabla u = (\partial_{x_i} u_j), 
$$
where $\lambda$ and $\mu$ are the Lam\'e coefficients, which satisfy
$$
\mu>0 \qquad \mbox{and} \qquad 2\mu + d\lambda \geq 0.
$$
\item Fourier's law $$q = \kappa \nabla \theta,$$ where $\kappa>0$ is the thermal conductivity.
\end{itemize}

We refer to~\cite{Lions} for a more detailed discussion of these assumptions. The equations become

\begin{equation}\notag \label{eq:cNS} \tag{cNS}
\begin{cases}
\partial_t \rho + {\rm div \, } \big( \rho u \big) 
= 0 , 
\\
\displaystyle 
\partial _t (\rho u) + 
{\rm div \,} \big( \rho u \otimes u \big) + \nabla ( \rho R \theta ) 
 = (\lambda + \mu) \nabla {\rm div \,} u + \mu \Delta u,
\\
\displaystyle 
\partial_t 
\Big[ \rho \Big( \frac{|u|^2}{2} + C_V \theta \Big) \Big] 
+ {\rm div} 
\Big[ u \Big( \rho \Big(\frac{|u|^2}{2} + C_V \theta \Big) 
      + \rho R \theta \Big) \Big]
- \kappa \Delta \theta
\\
\qquad \qquad \qquad \qquad \qquad \qquad \qquad \qquad \qquad= {\rm div \,} (\lambda ({\rm div \,} u) u + 2 \mu D(u) \cdot u ) .  
\end{cases}
\end{equation}

\subsection{Backward self-similar solutions} 
The equations \eqref{eq:cNS} exhibit a scaling invariance: 
The set of solutions is left invariant by the transformation 
\[
\rho(t,x) \to \rho (\lambda ^2 t , \lambda x) , 
\quad 
u(t,x) \to \lambda u (\lambda ^2 t, \lambda x) , \quad 
\quad 
\theta (t,x) \to \lambda ^2 \theta (\lambda ^2 t , \lambda x) , 
\quad \text{ for } \lambda > 0 .
\]
As a continuation of the previous paper~\cite{GI} which studies 
forward self-similar solutions in the radial case,  
we consider backward self-similar solutions, which has the following form. 
\begin{equation}\label{eq:back_selfsimilar}
\rho (t,x) = P \left( \frac{r}{\sqrt{T-t}} \right),  
\,\,
u(t,x) = \frac{1}{\sqrt{T-t}} U\left( \frac{r}{\sqrt {T-t}} \right)
  \frac{x}{r},
\,\, 
\theta (t,x) = \frac{1}{T-t} \Theta \left( \frac{r}{\sqrt {T-t}} \right), 
\end{equation}
where $T > 0$, $r = |x|$, $P, U , \Theta$ are scalar functions from 
$(0,\infty)$ to $\mathbb R$. 
It is natural to expect that there exist real numbers $P_\infty , U_\infty, 
\Theta_\infty$ such that 
\begin{equation}
\label{PUThetaInfinity}
P(r) \to P_\infty , \quad 
U(r) \sim \frac{U_\infty}{r} , \quad 
\Theta (r) \sim \frac{\Theta_\infty} {r^2}  
\quad \text{as } r \to \infty .
\end{equation}
This solution is then associated to self-similar blow up 
in $L^\infty$ of the velocity and the temperature, since 
\[
\lim _{t \to T}( u , \theta ) (t)
= \Big(  U_\infty \frac{x}{r^2} , \frac{\Theta}{r^2} \Big) . 
\]

\subsection{Known results} 

\subsubsection{Weak and strong solutions of~\eqref{eq:cNS}} 
In the very rich existing literature, we mention weak, finite-energy  solutions by Lions~\cite{Lions-1998}, variational solutions by 
Feireisl-Novotn\'y-Petzeltov\'a~\cite{FNP-2001} (see also Feireisl~\cite{F-2004}), classical solutions with finite energy by 
Matsumura-Nishida~\cite{MN-1980} (see also
Huang-Li~\cite{HL-2018} with vacuum), solutions in Besov spaces with the interpolation index one by Danchin~\cite{D-2001} (see also Chikami-Danchin~\cite{CD-2015}).

\subsubsection{Self-similar solutions of~\eqref{eq:cNS}}
There are only few results in this direction. 
Under a different scaling property from the parabolic type \eqref{eq:back_selfsimilar}, 
Qin-Su-Deng~\cite{QSD-2008} proved the non-existence of forward and backward self-similar solutions to the compressible Navier-Stokes equations in one dimension. Local energy of forward and backward self-similar solutions 
was also investigated in \cite{QSD-2008} 
but the total energy blows up 
at $t = 0$ and $t = T$, respectively, where $T$ is the given time 
appearing in  the definition of backward self-similar solutions. 
We also refer to related papers~\cite{GJ-2006} by Guo-Jiang 
(isothermal compressible Navier-Stokes equations) and 
Li-Chen-Xie~\cite{LCX-2013} (density-dependent viscosity).  
Finally, in a recent paper by the two first authors~\cite{GI}, the existence of forward
self-similar solutions has been established. 

\subsubsection{The case of incompressible Navier-Stokes} This case is different in two respects. First, the ansatz which we chose above 
(radial velocity) 
is incompatible with incompressibility, in fact, the velocity is irrotational; 
therefore, it is not possible to reduce the problem to a one-dimensional one, as we shall do in the present article. Second, the existence of forward self-similar solutions is known since strong solutions can be built up from  small self-similar data: see for instance Cannone and Planchon~\cite{CP}, Chemin~\cite{Chemin} and Koch and Tataru~\cite{KT}. The case of large self-similar data was recently treated by Jia and Sverak~\cite{JS}, who could prove the existence of smooth self-similar solutions.

It is known that the only possible backward self-similar solutions are the trivial ones
in several settings. The first result of this kind is due to Ne\v{c}as, R\r{u}\v{z}i\v{c}ka and \v{S}ver\'{a}k~\cite{NRS-1996} who proved that backward self-similar solutions with the velocity in $L^3 (\mathbb R^3)$ are zero. 
Tsai~\cite{Tsai-1998} generalized it under local energy bounds, and 
Chae and Wolf~\cite{CW-2017} extended the result in $L^{p,\infty}(\mathbb R^3)$ with $p > 3/2$.


\subsection{Analytical results} 
Consider solutions such that \eqref{eq:back_selfsimilar} is satisfied. 
The partial differential equations \eqref{eq:cNS} is equivalent to 
the following ordinary differential equations for 
$r = |x| > 0$: 
\begin{equation}\label{eq:ODEs_back}
\begin{cases}
\displaystyle 
\frac{1}{2} r P ' + P' U 
+ P \Big( U' + \frac{d-1}{r} U \Big) 
= 0 , 
\\[5mm]
\begin{split}
\displaystyle 
\frac{1}{2} 
 P  U 
& 
+ \frac{1}{2} r (P U)' 
+ (P U ^2) ' + \frac{d-1}{r} P U^2 
+(P R \Theta)'
\\
= 
& 
(2\mu + \lambda) 
 \Big( U'' + \frac{d-1}{r} U' - \frac{d-1}{r^2} U \Big) ,
\end{split}
\\[10mm]
\displaystyle 
\begin{split}
P 
 \Big( \frac{U^2}{2}
&  + C_V \Theta \Big) 
+ \frac{1}{2} r 
  \Big( P \Big( \frac{U^2}{2} + C_V \Theta \Big) \Big) '
+ \Big( U P \Big( \frac{U^2}{2} + C_V \Theta \Big) 
        + U P R \Theta \Big) '
\\
&
+ \frac{d-1}{r} 
  \Big( U P \Big( \frac{U^2}{2} + C_V \Theta \Big)
       + U P R \Theta
  \Big)
- \kappa \Big( \Theta '' + \frac{d-1}{r} \Theta ' \Big)
\\
= 
& 
2\mu \Big( (U')^2 + \frac{d-1}{r^2} U^2 \Big) 
+ \lambda \Big( U' + \frac{d-1}{r} U \Big) ^2
\\
& 
+ (2\mu + \lambda) 
 \Big( U'' + \frac{d-1}{r}U' - \frac{d-1}{r^2} U \Big) U. 
\end{split}
\end{cases}
\end{equation}

Hereafter, we always consider non-negative density $P \geq 0$ and 
temperature $\Theta \geq 0$ without explaining anything. 
We demonstrate a non-existence result of radially symmetric backward self-similar solutions 
under a smallness condition.

\begin{thm}\label{thm}
Let $d \geq 1$. 
There exists $\gamma > d$ such that the following hold. 
If 
there exist $0 < \beta < 1/2$ such that 
$U(r) \geq - \beta r$  for all $r > 0$,
\[
P(r)>0  \text{ almost every } r > 0 , \qquad 
{\liminf_{r\to \infty} P(r)r^{2-\varepsilon} > 0, } 
\]
for some $\varepsilon >0$, and 
\[
P,U,U',\Theta,\Theta' \in L^\infty (0,\infty), \quad 
P \in C^1 ((0,\infty)) , \quad 
U , \Theta \in C^2 ((0,\infty)), 
\]
\[
\sup_{r>0} 
\bigg( \frac{P\Theta}{ \min \{ 2\mu + \lambda, \kappa \} } 
   + \frac{Pr|U|}{\max \{ 2\mu+\lambda, \kappa \} } 
\bigg) 
\bigg( \gamma +  \frac{2\mu+\lambda}{\kappa}  
\bigg) ^{\log \gamma}
\ll 1, 
\]
then
$U \equiv \Theta \equiv 0, P = Constant$. 
\end{thm}
\noindent 
{\bf Remark.} 
\begin{itemize}

\item Notice that the smallness condition on the velocity profile $U$ is critical with respect to the scaling of the equation, while the smallness condition on the temperature profile is supercritical.

\item 
As seen from the proof of Theorem~\ref{thm}, we can relax the boundedness of $U, U'$ 
for large $r $ by polynomial growth of any order. 

\item 
The assumption $U(r) \geq - \beta r$ ($\beta < 1/2$) can be understood to justify 
the characteristic curve for the continuity equation. 

\item 
The smallness condition on $U$ is scaling invariant.
\end{itemize}

We can also obtain the following result for finite energy solutions.

\begin{thm}\label{thm:2}
Let $d = 1$ or $d \geq 3$. 
If 
\[
\Theta, P\Theta, P\Theta U , PU^2 \in L^1 (\mathbb R^d) , \quad U \in H^1 (\mathbb R^d) , 
\quad 
\lim_{\varepsilon \to 0} \varepsilon \int _{|x| \leq \varepsilon^{-1}} P|U|^3 ~dx = 0, 
\]
then $U \equiv \Theta \equiv 0$, $P \equiv Constant$.  
\end{thm}

\noindent 
{\bf Remark.} 
\begin{itemize}
\item  
Comparing the conditions in Theorem~\ref{thm:2} 
and \eqref{PUThetaInfinity}, the decay \eqref{PUThetaInfinity} 
implies finite energy only if $d = 1$.

\end{itemize}

\subsection{Numerical results} Section~\ref{LSNA} is devoted to a numerical study of the ODE~\eqref{eq:ODEs_back}. We follow a shooting method, by prescribing the values of $P$, $U$ and $\Theta$ at $\delta>0$, and solving for $r>\delta$ with the help of the Matlab solver ode45. Our ansatz for the values of $P$, $U$ and $\Theta$ is justified by the local existence result contained in~\cite{GI}, which ensures that we are (very close to) a smooth solution for $0<r<\delta$. We control the numerical stability of our approach by letting $\delta \to 0$.

The results of these simulations is very clear: for any set of values at $r=\delta$, the solution diverges very quickly - we refer to Section~\ref{LSNA} for more details.

This leads us to the 

\begin{conjecture} There are no backward self-similar solutions of~\eqref{eq:cNS} with radial symmetry. \end{conjecture}

One is tempted to ask: what if the radiality assumption is dropped?

\section{Proof of Theorem~\ref{thm}}

We use the functions
\begin{align*}
Z(r) & :=\frac{C_V}{\kappa} \int_0 ^r P(r_1)r_1 \bigg( \frac{1}{2} + \frac{U}{r_1} \bigg) dr_1, 
\qquad 
W(r)  := \frac{1}{2\mu+\lambda} \int_0 ^r P(r_1)r_1 \bigg( \frac{1}{2} + \frac{U}{r_1} \bigg) dr_1. 
\end{align*}
Let $m \in \mathbb N$ satisfy $2^{m+1} > d$. 
We divide into the two cases $\kappa \geq C_V (2\mu + \lambda)$ 
and $\kappa \leq C_V (2\mu + \lambda)$, and 
we show 3 kinds of energy inequalities in each case. 
We also use parameters $\delta _j > 0$ ($j = 1,2, \cdots , 6$) to be fixed later.
In this section, constants $C$ below are allowed to depend only on $R,C_V, d,m$ and of course their value changes from one line to another. 

\vskip3mm 

\noindent 
\subsection{The case when $\kappa \geq C_V (2\mu + \lambda)$}
\underline{Step 1}. Since the third equation of \eqref{eq:cNS} can be written
\begin{align*}
C_V \rho \big(\partial_t \theta + (u \cdot \nabla) \theta \big) + \rho R \theta \textrm{div} u - \kappa \Delta \theta = 2 \mu \vert D(u) \vert^2 + \lambda (\textrm{div} u)^2 ,
\end{align*}
we get with our ansatz
\begin{equation}\label{0731-1}
\begin{split}
& C_V P \Theta + \frac{C_V}{2} r P \Theta' + C_V PU \Theta' + PR \Theta \bigg(U' + \frac{d-1}{r} U \bigg) - \kappa \bigg(\Theta'' + \frac{d-1}{r} \Theta' \bigg) \\
& = 2 \mu \Bigg((U')^2 + \frac{d-1}{r^2} U^2 \Bigg) + \lambda \Bigg(U'+ \frac{d-1}{r} U \Bigg)^2.
\end{split}
\end{equation}
After multiplication by $e^{-Z(r)}$, it is written 
\begin{equation}\label{1002-1}
\begin{split}
& -\kappa \nabla \cdot \big(e^{-Z(r)} \nabla \Theta \big) 
+ C_V P\Theta e^{-Z(r)} + PR\Theta \Big( U' + \frac{d-1}{r} U \Big) e^{-Z(r)} \\
 &=2 \mu \bigg( (U')^2 + \frac{d-1}{r^2} U^2 \bigg) e^{-Z(r)} + \lambda e^{-Z(r)} \bigg( U' + \frac{d-1}{r} U \bigg)^2.
\end{split}
\end{equation}
We multiply $\Theta^{2^m -1} $, integrate over the domain, and integrate by parts in the left-hand side 
and get that 
\begin{align*}
& 
\int_{\mathbb R^d} 
\bigg\{ 
(2^m -1)\kappa  |\nabla \Theta|^2 \Theta ^{2^{m}-2}
+ C_V P \Theta ^{2^m} + PR\Theta ^{2^m}  \Big( U' + \frac{d-1}{r}{U} \Big) 
\bigg\} e^{-Z(r)} dx 
\\
& 
= \int_{\mathbb R^d} 
\bigg\{ 2\mu \Big( (U')^2 + \frac{d-1}{r^2}U^2 \Big) \Theta ^{2^m-1}
    + \lambda \Big(U' + \frac{d-1}{r}U \Big)^2 \Theta^{2^m -1} \bigg\} e^{-Z(r)} dx .
\end{align*}
We have 
\begin{align*}
 PR\Theta ^{2^m}  \Big( U' + \frac{d-1}{r}{U} \Big) 
& \leq 
\frac{R^2 P \Theta}{2\mu + \lambda} P \Theta ^{2^m}
+ \frac{2\mu+\lambda}{4} \Theta^{2^{m} -1} \Big( U' + \frac{d-1}{r}{U} \Big) ^2 ,
\end{align*}
which implies the first energy inequality 
\begin{equation}\label{0730-1}
\begin{split}
& 
\int_{\mathbb R^d} 
\bigg\{ 
(2^m -1)\kappa  |\nabla \Theta|^2 \Theta^{2^m -2} 
+ \Big( C_V - \frac{R^2 P \Theta}{2\mu + \lambda} \Big) P \Theta ^{2^m} 
\bigg\} e^{-Z(r)} dx 
\\
& 
\leq \int_{\mathbb R^d} 
\bigg\{ { \mu(d+2)} \Big( (U')^2 + \frac{d-1}{r^2}U^2 \Big) \Theta ^{2^m -1}
    + 2\lambda \Big(U' + \frac{d-1}{r}U \Big)^2 \Theta ^{2^m -1} \bigg\} e^{-Z(r)} dx .
\end{split}
\end{equation}

\bigskip

\noindent \underline{Step 2.} 
Since the second equation of \eqref{eq:ODEs_back} can be written
\begin{equation}\label{1003-1}
\frac{1}{2} PU + P\Big( \frac{1}{2}r + U \Big) U' + (PR\Theta)' 
 - (2\mu + \lambda) \Big( U'' + \frac{d-1}{r}U' - \frac{d-1}{r^2} U \Big) = 0
\end{equation}
or
\begin{equation}\label{0730-2}
\Big( \frac{1}{2} PU + (PR\Theta)' \Big) e^{-W(r)} 
 + (2\mu + \lambda) \left\{ 
   -\nabla \cdot  \big( e^{-W(r)} \nabla U \Big) + \frac{d-1}{r^2} U e^{-W(r)} \right\}
   = 0.
\end{equation}
For $n = 1,2, \cdots , m$, after multiplication by $\Theta^{2^m - 2^{n-1} } U^{2^n -1} e^{W(r)-Z(r)}$, 
integration over $\mathbb R^d$, integration by parts, 
it follows that 
\[
\begin{split}
& 
\int_{\mathbb R^d} 
\Big( \frac{1}{2}PU^{2^n} \Theta^{2^m- 2^{n-1} } + (PR\Theta)' U^{2^n -1} \Theta^{2^m- 2^{n-1} } \Big) e^{-Z(r)} dx
\\
& 
+ (2\mu + \lambda) \int_{\mathbb R^d} 
  \Big( (2^n -1)|\nabla U|^2 U^{2^n -2}\Theta^{2^m - 2^{n-1}} 
      + (2^m - 2^{n-1}) U^{2^n -1} (\nabla U \cdot \nabla \Theta) \Theta ^{2^m - 2^{n-1} -1} 
      \\
& \qquad \qquad \qquad    + \frac{d-1}{r^2}U^{2^n} \Theta^{2^m- 2^{n-1} } 
  \Big) e^{-Z(r)}dx
\\
& 
+ \Big( 1 - \frac{C_V (2\mu+\lambda)}{\kappa} \Big) \int_{\mathbb R^d} 
  U' U^{2^n -1} \Theta^{2^m- 2^{n-1} } P\Big( \frac{r}{2} + U \Big) e^{-Z(r)}dx
  = 0.
\end{split}
\]
We estimate the pressure term and the cross term 
\begin{align*}
| P 'R\Theta U^{2^n -1} \Theta^{2^m- 2^{n-1} } | 
& 
= P \frac{U'+ \frac{d-1}{r}U}{\frac{1}{2}r + U} R \Theta U^{2^n -1} \Theta^{2^m- 2^{n-1} }
\\
& \leq  \, C (P\Theta) U^{2^n -2} \Theta^{2^m- 2^{n-1} } \Big( (U')^2 + \frac{d-1}{r^2}U^2 \Big) ,
\\
| P R\Theta' U^{2^n -1} \Theta^{2^m- 2^{n-1} } | 
& \leq  C \delta_1 \kappa  (\Theta')^2 \Theta^{2^m -2} 
     +  \frac{1}{\delta _1 \kappa} P ^2 U^{2^{n+1}-2} \Theta^{2^m - 2^n +2} 
\\
& \leq C \delta_1 \kappa  (\Theta')^2 \Theta^{2^m -2} 
     +  \frac{P\Theta }{\delta _1 \kappa} 
      \Big( PU^{2^n} \Theta ^{2^m - 2^{n-1}} + PU^{2^{n+1}} \Theta ^{2^m - 2^{n}} \Big) , 
\quad  (\delta_1 > 0)
\end{align*}
\[
\begin{split}
(2\mu + \lambda)
& 
(2^m - 2^{n-1}) 
U^{2^n -1} (\nabla U \cdot \nabla \Theta) \Theta ^{2^m - 2^{n-1} -1} 
\\
\leq 
& (2^m - 2^{n-1})\frac{\delta _2 \kappa }{2} |\nabla \Theta |^2 \Theta ^{2^m -2}
+ (2^m - 2^{n-1}) \frac{(2\mu+\lambda)^2 }{2 \delta _2 \kappa}  
    |\nabla U| ^2 U^{2^{n+1} -2} \Theta ^{2^m - 2^n} .
\quad (\delta _2 > 0)
\end{split}
\]
and write by the integration by parts 
\begin{align*}
& 
\Big( 1 - \frac{C_V (2\mu+\lambda)}{\kappa} \Big) \int_{\mathbb R^d} 
  U' U^{2^n-1} \Theta^{2^{m} - 2^{n-1}} P\Big( \frac{r}{2} + U \Big) e^{-Z(r)}dx
\\
& = 
\Big( 1 - \frac{C_V (2\mu+\lambda)}{\kappa} \Big) \int_{\mathbb R^d} 
\bigg\{ \frac{C_V}{\kappa}\frac{U^{2^n}}{2^n} \Theta^{2^m - 2^{n-1}} P^2 \Big( \frac{r}{2} + U \Big) ^2 
   - \frac{d}{2^{n+1}} PU^{2^n} \Theta^{2^m - 2^{n-1}} 
\\
& \qquad 
   - (2^m - 2^{n-1}) \frac{U^{2^n}}{2^n} \Theta' \Theta ^{2^m - 2^{n-1} -1} P \Big( \frac{r}{2}+U \Big) 
\bigg\} e^{-Z(r)} dx . 
\end{align*}
By applying the Schwarz inequality to the last term above, 
\begin{align*}
& 
\bigg| 
\Big( 1 - \frac{C_V (2\mu+\lambda)}{\kappa} \Big) \int_{\mathbb R^d} 
(2^m - 2^{n-1}) \frac{U^{2^n}}{2^n} \Theta' \Theta ^{2^m - 2^{n-1}-1} P \Big( \frac{r}{2}+U \Big) 
e^{-Z(r)} \bigg| dx
\\
& 
\leq 
(2^m - 2^{n-1})\int_{\mathbb R^d} \frac{\delta _3}{2} \kappa |\Theta'|^2 \Theta ^{2^m -2} e^{-Z(r)} 
\\
& \quad 
+ \frac{2^m - 2^{n-1}}{2^n}  \frac{1}{\kappa \delta _3 }
\Big( 1 - \frac{C_V (2\mu+\lambda)}{\kappa} \Big) \int_{\mathbb R^d} 
  \frac{U^{2^{n+1}} }{2^{n+1}}
  \Theta ^{2^m - 2^n} P^2 \Big( \frac{r}{2}+U\Big)^2 e^{-Z(r)} dx. 
 \qquad (\delta _3 > 0) 
\end{align*}
and furthermore, 
\[
\begin{split}
& 
PU^{2^n} \Theta ^{2^m - 2^{n-1}}
\leq \frac{1}{\delta _4} PU^{2^{m+1}} + \delta _4 P\Theta ^{2^m} , 
\qquad (\delta _4 > 0)
\\
& 
PU^{2^{n+1}} \Theta ^{2^m - 2^{n}}
\leq \frac{1}{\delta _4} PU^{2^{m+1}} + \delta _4 P\Theta ^{2^m} . 
\end{split}
\]
We then obtain the second energy inequality 
\begin{equation}\label{0730-3}
\begin{split}
& 
 (2\mu + \lambda) \int_{\mathbb R^d} 
  \Big\{ \Big(2^n -1 - \frac{CP\Theta}{2\mu + \lambda}\Big)|\nabla U|^2 U^{2^n -2}\Theta^{2^m - 2^{n-1}} 
      \\
& \qquad \qquad \qquad    
  + \Big( 1 - \frac{CP\Theta}{2\mu+\lambda} \Big) \frac{d-1}{r^2}U^{2^n} \Theta^{2^m- 2^{n-1} } 
  \Big\} e^{-Z(r)}dx
\\
& 
+ \Big( 1 - \frac{C_V (2\mu+\lambda)}{\kappa} \Big) \int_{\mathbb R^d} 
  \frac{C_V}{\kappa}\frac{U^{2^n}}{2^n} \Theta^{2^m - 2^{n-1}} P^2 \Big( \frac{r}{2} + U \Big) ^2 
  e^{-Z(r)}dx
\\
\leq 
& \int_{\mathbb R^d} 
  \bigg[ 
 \Big\{ C \delta _1 + (2^m - 2^{n-1})\frac{\delta _2 +\delta_3}{2} \Big\} \kappa |\Theta ' |^2 \Theta ^{2^m -2}
\\
& + (2^m - 2^{n-1}) \frac{(2\mu+\lambda)^2}{2 \delta _2 \kappa}  
    |\nabla U| ^2 U^{2^{n+1} -2} \Theta ^{2^m - 2^n}
\\
& + \frac{2^m - 2^{n-1}}{2^n}  \frac{1}{\kappa \delta _3 }
\Big( 1 - \frac{C_V (2\mu+\lambda)}{\kappa} \Big) 
  \frac{U^{2^{n+1}} }{2^{n+1}}
  \Theta ^{2^m - 2^n} P^2 \Big( \frac{r}{2}+U\Big)^2 
\\
& + \Big\{ \frac{P\Theta}{\delta_1 \kappa} + 
    \Big| \frac{1}{2} - \frac{P\Theta}{\delta _1 \kappa} 
        - \frac{d}{2^{n+1}} \Big( 1- \frac{C_V(2\mu+\lambda)}{\kappa} \Big) 
    \Big| 
    \Big\}  
    \Big( \frac{1}{\delta _4} PU^{2^{m+1}} + \delta _4 P\Theta ^{2^m}\Big)
  \bigg] e^{-Z(r)}dx . 
\end{split}
\end{equation}

\noindent \underline{Step 3}. 
We multipy \eqref{0730-2} by $U^{2^{m+1}-1} e^{W(r) - Z(r)}$, integrate over the domain 
and then obtain 
\[
\begin{split}
& 
\int_{\mathbb R^d} 
\bigg( 
\frac{1}{2} PU^{2^{m+1}} +(PR\Theta)' U^{2^{m+1}-1} 
\\
& + (2\mu +\lambda)\Big( (2^{m+1}-1) |\nabla U|^2 U^{2^{m+1}-2} + \frac{d-1}{r^2}U^{2^{m+1}} \Big) 
\bigg) e^{-Z(r)}dx
\\
& 
+ 
\Big( 1- \frac{C_V(2\mu+\lambda)}{\kappa} \Big) \int _{\mathbb R^d} 
U' U^{2^{m+1}-1} P\Big( \frac{r}{2}+U \Big) e^{-Z(r)}dx
= 0 .
\end{split}
\]
On the pressure term, it follows from the integration by parts that 
\[
\begin{split}
& 
 \int_{\mathbb R^d} (PR\Theta)' U^{2^{m+1} -1}  e^{-Z(r)} 
\\
= 
& 
\int_{\mathbb R^d} 
 \Big\{ - PR\Theta (2^{m+1} -1)U^{2^{m+1} -2} U' 
      + PR\Theta U^{2^{m+1}-1} \frac{C_V}{\kappa} P \Big( \frac{r}{2}+U\Big) 
\\
&      - PR\Theta U^{2^{m+1} -1 } \frac{d-1}{r}
 \Big\}e^{-Z(r)}dx
\end{split}
\]
and we estimate 
\[
\begin{split}
| PR\Theta (2^{m+1} -1)U^{2^{m+1} -2} U' | 
\leq 
&\frac{C}{\delta_5 (2\mu+\lambda)} P^2 \Theta ^2 U^{2^{m+1}-2} 
    + (2^{m+1} -1) \delta_5 (2\mu +\lambda ) |U'|^2 U^{2^{m+1}-2}
\\
\leq 
&\frac{ C P \Theta \Big( P\Theta ^{2^{m}} + P U ^{2^{m+1}} \Big) }{\delta_5 (2\mu+\lambda)} 
    + (2^{m+1} -1) \delta_5 (2\mu +\lambda ) |U'|^2 U^{2^{m+1}-2} , 
\\
& (\delta _5 > 0)
\\
\Big| PR\Theta U^{2^{m+1}-1} \frac{C_V}{\kappa} P \Big( \frac{r}{2}+U\Big) 
\Big| 
\leq 
& \frac{C_VR PrU}{4\kappa} (P\Theta ^{2^m} + P U^{2^{m+1}}) 
  + \frac{C_V R P\Theta}{\kappa} P U^{2^{m+1}}
\\
PR\Theta U^{2^{m+1} -1 } \frac{d-1}{r}
\leq 
& \frac{C P \Theta}{\delta_6 (2\mu+\lambda)} 
   \Big( P\Theta ^{2m} + P U^{2^{m+1}} \Big)  
   + \delta_6 (2\mu + \lambda) \frac{d-1}{r}U^{2^{m+1}} .
\\
& (\delta _6 > 0)
\end{split}
\]
It follows from the integration by parts 
that 
\[
\begin{split}
& 
\Big( 1- \frac{C_V(2\mu+\lambda)}{\kappa} \Big) \int _{\mathbb R^d} 
U' U^{2^{m+1} -1} P\Big( \frac{r}{2}+U \Big) e^{-Z(r)}
\\
= 
& 
\Big( 1- \frac{C_V(2\mu+\lambda)}{\kappa} \Big) \frac{C_V}{\kappa} \int _{\mathbb R^d} 
\frac{U^{2^{m+1}}}{2^{m+1}} P^2 \Big( \frac{r}{2}+U \Big) ^2 e^{-Z(r)} 
\\
& - \Big( 1- \frac{C_V(2\mu+\lambda)}{\kappa} \Big) \frac{d}{2^{m+2}} \int _{\mathbb R^d} 
PU^{2^{m+1}} e^{-Z(r)} .
\end{split}
\]
Therefore, we obtain the third energy inequality 
\begin{equation}\label{0730-4}
\begin{split}
& 
\int_{\mathbb R^d} 
\bigg[ 
\Big\{ \frac{1}{2} - \frac{C P \Theta}{\delta_5 (2\mu+\lambda)} 
   - \frac{C_VR PrU}{4\kappa} -\frac{C_V R P\Theta}{\kappa} 
   - \frac{C P \Theta}{\delta_6 (2\mu+\lambda)} 
\\
& \qquad \qquad    - \Big( 1- \frac{C_V(2\mu+\lambda)}{\kappa} \Big)  \frac{d}{2^{m+2}}
\Big\} PU^{2^{m+1}}
\\
& + (2\mu +\lambda)
    \Big( (2^{m+1}-1) (1 - \delta _5) |\nabla U|^2 U^{2^{m+1}-2} 
         + (1-\delta _6 )\frac{d-1}{r^2}U^{2^{m+1}} 
    \Big) 
\bigg] e^{-Z(r)}
\\
& 
+ 
\Big( 1- \frac{C_V(2\mu+\lambda)}{\kappa} \Big) \frac{C_V}{\kappa} \int _{\mathbb R^d} 
\frac{U^{2^{m+1}}}{2^{m+1}} P^2 \Big( \frac{r}{2}+U \Big) ^2 e^{-Z(r)} 
\\
\leq 
& \int_{\mathbb R^d} 
\bigg[ \Big( \frac{C P \Theta}{\delta_5 (2\mu+\lambda)} + \frac{C_VR PrU}{4\kappa} 
             + \frac{C P \Theta}{\delta_6 (2\mu+\lambda)}  \Big) 
     P\Theta ^{2^m} 
\bigg]e^{-Z(r)},
\end{split}
\end{equation}
where $\delta_5=\delta_6=\frac{1}{2}.$ \\

\bigskip

\noindent \underline{Step 4}.
We denote
\begin{align*}
I^{(m,n)} &:= \int_{\mathbb{R}^d} \vert \nabla U \vert^2 U^{2^n-2} \Theta^{2^m-2^{n-1}} e^{-Z(r)} dx,~1\leqslant n \leqslant m+1, \\
J^{(m,n)} &:= \int_{\mathbb{R}^d} \frac{d-1}{r^2} U^{2^n} \Theta^{2^m-2^{n-1}} e^{-Z(r)} dx,~1\leqslant n \leqslant m+1, \\
K^{(m,n)} &:= \int_{\mathbb{R}^d} U^{2^n} \Theta^{2^m-2^{n-1}} P^2 \bigg(\frac{r}{2}+U \bigg)^2 e^{-Z(r)} dx,~1\leqslant n \leqslant m+1, \\
L^{(m)} &:= \int_{\mathbb{R}^d} P U^{2^{m+1}} e^{-Z(r)} dx, \\
M^{(m)} &:= \int_{\mathbb{R}^d} P \Theta^{2^{m}} e^{-Z(r)} dx , \\
N^{(m)} &:= \int_{\mathbb{R}^d} \vert \nabla \Theta \vert^2 \Theta^{2^m-2} e^{-Z(r)} dx .
\end{align*}
Let $A, B, D>0$ to be determined later. Consider the sum
\begin{align*}
\eqref{0730-1} 
+ \sum_{n=1}^{m} B^{n} A^{2^{n-1}-1} \times \eqref{0730-3}_n
+ D \times \eqref{0730-4},
\end{align*}
where we choose
\begin{align*}
\delta_1 = \delta_2 = \delta_3 = \delta_4 = \frac{1}{A^{2^{n-1}}}
\end{align*}
in $\eqref{0730-3}_n.$
\\
This yields
\begin{align} \label{ineq-concl}
c_1 L^{(m)} + c_2 M^{(m)} + c_3 N^{(m)} + \sum_{n=1}^m \Bigg( c_{4,n} I^{(m,n)} + c_{5,n} J^{(m,n)} + c_{6,n} K^{(m,n)} \Bigg) \leqslant 0,
\end{align}
where the constants are defined as:
\begin{align*}
c_1 & := D \Bigg \lbrace \frac{1}{2} - \bigg( \frac{2C}{2\mu+\lambda} + \frac{C_V R}{\kappa} + \frac{2C}{2\mu+\lambda} \bigg) \sup_{r>0} P \Theta -  \frac{C_V R}{4 \kappa} \sup_{r>0} P r \vert U \vert - \bigg(1- \frac{C_V (2\mu+\lambda)}{\kappa} \bigg) \frac{d}{2^{m+2}} \Bigg \rbrace \\
& - \Bigg \lbrace \sum_{n=1} ^m \frac{B^{n} A^{3 \times 2^{n-1}-1}}{\kappa}  \sup_{r>0} P \Theta + B^{n} A^{2^n -1} \sup_{r>0} \bigg \vert  \frac{1}{2} - \frac{P \Theta}{\kappa} A^{2^{n-1}} - \frac{d}{2^{n-1}} \bigg(1- \frac{C_V (2\mu+\lambda)}{\kappa} \bigg) \bigg \vert \Bigg \rbrace , 
\end{align*}
\begin{align*}
c_2 & := \bigg(C_V - \frac{R^2}{2\mu+\lambda} \sup_{r>0} P \Theta \bigg) - D \Bigg[ \frac{4C}{2\mu+\lambda}  \sup_{r>0} P \Theta + \frac{C_V R}{4 \kappa} \sup_{r>0} P r \vert U \vert \Bigg]  \\
 & - \frac{1}{A} \sum_{n=1} ^m B^{n} \Bigg( \frac{A^{2^{n-1}}}{\kappa} \sup_{r>0} P \Theta + \sup_{r>0} \bigg \vert  \frac{1}{2} - \frac{P \Theta}{\kappa} A^{2^{n-1}} - \frac{d}{2^{n-1}} \bigg(1- \frac{C_V (2\mu+\lambda)}{\kappa} \bigg) \bigg \vert \Bigg),
\end{align*}
\begin{align*}
c_3 & := (2^m -1) \kappa - \frac{1}{A} \sum_{n=1}^m B^{n} \big( C + 2^m-2^{n-1} \big) \kappa  .
\end{align*}
For the constants $c_{4,n}$ there are different cases to consider:
\begin{align*}
c_{4,1} & := (2 \mu + \lambda) B \big(1 - \frac{C}{2\mu+\lambda} \sup_{r>0} P \Theta \big) -(d+2)\mu - 4 \lambda , 
\end{align*}
when $1< n < m+1,$
\begin{align*}
c_{4,n} & :=A^{2^{n-1}-1} B^{n-1} \Bigg( (2\mu+\lambda) B \big(2^n -1 -  \frac{C}{2\mu+\lambda} \sup_{r>0} P \Theta  \big) - (2^m - 2^{n-2}) \frac{(2\mu+\lambda)^2}{2 \kappa} \Bigg) , 
\end{align*}
and finally
\begin{align*}
c_{4,m+1} & := D \frac{2 \mu+\lambda}{2} (2^{m+1}-1) - 2^{m-1} \frac{(2\mu+\lambda)^2}{2 \kappa} B^{m} A^{2^{m}-1}  
\end{align*}
Similarly, we have different cases for $c_{5,n}:$
\begin{align*}
c_{5,1} :=  (2\mu+\lambda) B \big( 1 - \frac{C}{2\mu+\lambda} \sup_{r>0} P \Theta \big) -(d+2) \mu - 4 \lambda (d-1)       ,
\end{align*}
for $1<n<m+1$ we have the formula
\begin{align*}
c_{5,n} := (2\mu+\lambda) B^{n} A^{2^{n-1} -1}  \big( 1 - \frac{C}{2\mu+\lambda} \sup_{r>0} P \Theta \big) ,
\end{align*}
and finally 
\begin{align*}
c_{5,m+1} :=  D \frac{2\mu+\lambda}{2}.
\end{align*}
We conclude with the last set of constants:
\begin{align*}
c_{6,1} :=  B \big( 1 - \frac{C_V(2\mu+\lambda)}{\kappa} \big) \frac{C_V}{2\kappa}    ,
\end{align*}
for $1 < n < m+1$ we have
\begin{align*}
c_{6,n} :=  B^{n-1} A^{2^{n-1}-1} \bigg( B \big( 1 - \frac{C_V(2\mu+\lambda)}{\kappa} \big) \frac{C_V}{2^n \kappa} - \frac{2^m - 2^{n-2}}{2^{n-1}} \frac{1}{2^n \kappa} \big(1- \frac{C_V(2\mu+\lambda)}{\kappa} \big)  \bigg)  ,
\end{align*}
and finally
\begin{align*}
c_{6,m+1}:=  D  \big(1-\frac{C_V(2\mu+\lambda)}{\kappa} \big) \frac{C_V}{2^{m+1}\kappa} - \frac{1}{2^{m+2} \kappa} \big(1 - \frac{C_V(2\mu+\lambda)}{\kappa}  \big) B^m A^{2^{m}-1}    .
\end{align*}
We are going to conclude the proof by showing that for an appropriate choice of $A, B, D$ and our smallness assumptions, all the constants appearing in \eqref{ineq-concl} are strictly positive. \\
First note that without any condition, $c_{5,m+1} , c_{6,1} >0.$ \\
The condition $\sup_{r>0} P \Theta \ll 1 $ (here the implicit constants depend on the physical constants) guarantees that $ c_{5,n}>0 $ for $1<n<m+1.$ \\
Taking in addition $1 \ll B$ allows us to conclude that $c_{5,1} >0, c_{4,n} >0$ for $1 \leqslant n \leqslant m$ and $c_{6,n}>0$ for $2 \leqslant n \leqslant m.$ \\
Fixing $B^m \ll A$ forces $c_3 >0.$ Requiring $B^m A^{3 \times 2^m} \ll D$ gives $c_1, c_{4,m+1}, c_{6,m+1} >0.$
\\
The only constant remaining is $c_2,$ but the smallness of $\sup_{r>0} \big( P \Theta  + P r \vert U \vert \big)$ combined with the choice of $A$ above guarantees that it is strictly positive. \\
We conclude that $U = \Theta = 0$ and $P = Constant.$

\subsection{ The case when $\kappa \leq C_V (2\mu + \lambda)$}.
We consider a similar argument to the previous case with the weight $e^{-W(r)}$. 

We multiply the equation \eqref{1002-1} 
by $\Theta^{2^m -1} e^{Z(r) - W(r)}$, integrate over the domain, and integrate by parts in the left hand side 
and get that 
\begin{align*}
& 
\int_{\mathbb R^d} 
\bigg\{ 
(2^m -1)\kappa  |\nabla \Theta|^2 \Theta ^{2^{m}-2}
+ C_V P \Theta ^{2^m} + PR\Theta ^{2^m}  \Big( U' + \frac{d-1}{r}{U} \Big) 
\bigg\} e^{-W(r)} dx 
\\
& + C_V \Big(1 - \frac{\kappa}{C_V(2\mu+\lambda)} \Big) 
   \int _{\mathbb R^d}  \Theta' \Theta ^{2^m -1} P\Big( \frac{r}{2} + U \Big) e^{-W(r)} dx 
\\
& 
= \int_{\mathbb R^d} 
\bigg\{ 2\mu \Big( (U')^2 + \frac{d-1}{r^2}U^2 \Big)^2 \Theta ^{2^m-1}
    + \lambda \Big(U' + \frac{d-1}{r}U \Big)^2 \Theta^{2^m -1} \bigg\} e^{-W(r)} dx .
\end{align*}
We have from the Schwarz inequality that 
\begin{align*}
 PR\Theta ^{2^m}  \Big( U' + \frac{d-1}{r}{U} \Big) 
& \leq 
\frac{R^2 P \Theta}{2\mu + \lambda} P \Theta ^{2^m}
+ \frac{2\mu+\lambda}{4} \Theta^{2^{m} -1} \Big( U' + \frac{d-1}{r}{U} \Big) ^2 ,
\end{align*}
and the second integral in the left-hand side is equal to
\[
 C_V \Big(1 - \frac{\kappa}{C_V(2\mu+\lambda)} \Big) 
   \int _{\mathbb R^d}  
   \Big\{ \frac{1}{2\mu + \lambda}\frac{\Theta ^{2^m}}{2^m} P^2 \Big( \frac{r}{2} + U \Big) ^2 
       - \frac{d}{2^{m+1}} \Theta ^{2^m} 
   \Big\} e^{-W(r)} ~dx 
\]
which imply the first energy inequality 
\begin{equation}\label{0730-5}
\begin{split}
& 
\int_{\mathbb R^d} 
\bigg[ 
(2^m -1)\kappa  |\nabla \Theta|^2 \Theta^{2^m -2} 
\\
& + \Big\{ C_V - \frac{R^2 P \Theta}{2\mu + \lambda} 
  - C_V \Big( 1- \frac{\kappa}{C_V (2\mu+\lambda)} \Big) \frac{d}{2^{m+1}}  \Big\} P \Theta ^{2^m} 
\\
& + C_V \Big(1 - \frac{\kappa}{C_V(2\mu+\lambda)} \Big) 
    \frac{1}{2\mu + \lambda}\frac{\Theta ^{2^m}}{2^m} P^2 \Big( \frac{r}{2} + U \Big) ^2 
   \bigg] e^{-W(r)} dx 
\\
& 
\leq \int_{\mathbb R^d} 
\bigg\{ 3\mu \Big( (U')^2 + \frac{d-1}{r^2}U^2 \Big) \Theta ^{2^m -1}
    + 2\lambda \Big(U' + \frac{d-1}{r}U \Big)^2 \Theta ^{2^m -1} \bigg\} e^{-W(r)} dx .
\end{split}
\end{equation}
To obtain the second energy inequality, 
we multiply \eqref{0730-2} by $\Theta^{2^m - 2^{n-1} } U^{2^n -1} e^{-W(r)}$, 
estimate analogously to the proof of \eqref{0730-3} and then obtain 
\begin{equation}\label{0730-6}
\begin{split}
& 
 (2\mu + \lambda) \int_{\mathbb R^d} 
  \Big\{ \Big(2^n -1 - \frac{CP\Theta}{2\mu + \lambda}\Big)|\nabla U|^2 U^{2^n -2}\Theta^{2^m - 2^{n-1}} 
      \\
& \qquad \qquad \qquad    
  + \Big( 1 - \frac{CP\Theta}{2\mu+\lambda} \Big) \frac{d-1}{r^2}U^{2^n} \Theta^{2^m- 2^{n-1} } 
  \Big\} e^{-W(r)}dx
\\
\leq 
& \int_{\mathbb R^d} 
  \bigg[ 
 \Big\{ C \delta _1 + (2^m - 2^{n-1})\frac{\delta _2 }{2} \Big\} \kappa |\Theta ' |^2 \Theta ^{2^m -2}
\\
& + (2^m - 2^{n-1}) \frac{(2\mu+\lambda)^2}{2 \delta _2 \kappa}  
    |\nabla U| ^2 U^{2^{n+1} -2} \Theta ^{2^m - 2^n}
\\
& + \Big\{ \frac{P\Theta}{\delta_1 \kappa} + 
    \Big| \frac{1}{2} - \frac{P\Theta}{\delta _1 \kappa} 
    \Big| 
    \Big\}  
    \Big( \frac{1}{\delta _4} PU^{2^{m+1}} + \delta _4 P\Theta ^{2^m}\Big)
  \bigg] e^{-W(r)} dx. 
\end{split}
\end{equation}
For the third energy inequality, 
we multipy \eqref{0730-2} by $U^{2^{m+1}-1} e^{-W(r)} $, 
and analogously to the proof of \eqref{0730-4}, 
\begin{equation}\label{0730-7}
\begin{split}
& 
\int_{\mathbb R^d} 
\bigg[ 
\Big\{ \frac{1}{2} - \frac{C P \Theta}{\delta_5 (2\mu+\lambda)} 
   - \frac{R PrU}{4(2\mu+\lambda)} -\frac{R P\Theta}{2\mu+\lambda} 
   - \frac{C P \Theta}{\delta_6 (2\mu+\lambda)} 
\Big\} PU^{2^{m+1}}
\\
& + (2\mu +\lambda)
    \Big( (2^{m+1}-1) (1 - \delta _5) |\nabla U|^2 U^{2^{m+1}-2} 
         + (1-\delta _6 )\frac{d-1}{r^2}U^{2^{m+1}} 
    \Big) 
\bigg] e^{-W(r)}dx
\\
\leq 
& \int_{\mathbb R^d} 
\bigg[ \Big( \frac{C P \Theta}{\delta_5 (2\mu+\lambda)} + \frac{R PrU}{4(2\mu + \lambda)} 
             + \frac{C P \Theta}{\delta_6 (2\mu+\lambda)}  \Big) 
     P\Theta ^{2^m} 
\bigg]e^{-W(r)} dx.
\end{split}
\end{equation}
Exactly as we did in the previous case, we consider the following combination of inequalities
\[
\eqref{0730-5} 
+ \sum_{n=1}^{m} B^n A^{2^{n-1}-1} \times \eqref{0730-6}_n 
+ D \times \eqref{0730-7} ,
\]
for $A, B ,D >0$ appropriately chosen. This leads to $U \equiv \Theta \equiv 0$, and $P = Constant$.

\section{Proof of Theorem~\ref{thm:2}}
We rewrite the equation \eqref{0731-1} into 
\begin{equation}\notag 
\begin{split}
& C_V P \Theta + P \Big( \frac{r}{2} + U \Big) (C_V \Theta) ' 
+ PR \Theta {\rm \, div \, } \bigg(U \frac{x}{r} \bigg) 
- \kappa \Delta \Theta 
\\
& = 2 \mu \Big((U')^2 + \frac{d-1}{r^2} U^2 \Big) + \lambda \Big(U'+ \frac{d-1}{r} U \Big)^2.
\end{split}
\end{equation}
We note 
\[
-(\Delta U) U = -\Delta \frac{U^2}{2} + |\nabla U |^2
\]
and multiply the momentum equation \eqref{1003-1} by $U$ 
and then get 
\[
\begin{split}
\displaystyle 
& P  \frac{U^2}{2}
+ P \Big( \frac{r}{2} + U \Big) \Big( \frac{U^2}{2} \Big) ' 
+ \Big(\nabla (P R \Theta)\Big) \cdot U \frac{x}{r}
\\
& 
- (2\mu + \lambda) \Delta \frac{U^2}{2} 
+ (2\mu + \lambda) \Big( |\nabla U|^2 +  \frac{d-1}{r^2} U^2\Big) 
= 0. 
\end{split}
\]
The sum of the two equations is 
\[
\begin{split}
& 
P \Big( C_V \Theta + \frac{U^2}{2} \Big) 
+ P\Big( \frac{r}{2} + U \Big) \Big( C_V \Theta + \frac{U^2}{2} \Big) ' 
+ {\rm \, div \,} \Big( PR\Theta U \frac{x}{r} \Big) 
\\
& 
- \kappa \Delta \Theta 
- (2\mu + \lambda) \Delta \frac{U^2}{2} 
+ (2\mu + \lambda) \Big( |\nabla U|^2 +  \frac{d-1}{r^2} U^2\Big) 
\\
= 
&  
 2 \mu \Big((U')^2 + \frac{d-1}{r^2} U^2 \Big) + \lambda \Big(U'+ \frac{d-1}{r} U \Big)^2. 
\end{split}
\]
The right-hand side is written 
\[
\begin{split}
& 
 2 \mu \Big((U')^2 + \frac{d-1}{r^2} U^2 \Big) + \lambda \Big(U'+ \frac{d-1}{r} U \Big)^2
\\
 =
 & (2\mu+\lambda) \Big( |\nabla U|^2 + \frac{d-1}{r^2} U^2 \Big)
  + \lambda \Big( \frac{d-1}{r} (U^2) ' + \frac{(d-1)(d-2)}{r^2} U^2 \Big) ,
\\
 =
 & (2\mu+\lambda) \Big( |\nabla U|^2 + \frac{d-1}{r^2} U^2 \Big)
  + \lambda {\rm \, div \,} \Big( \frac{d-1}{r} U^2 \frac{x}{r} \Big) ,
\end{split}
\]
which implies 
\begin{equation*}
\begin{split}
& 
P \Big( C_V \Theta + \frac{U^2}{2} \Big) 
+ P\Big( \frac{r}{2} + U \Big) \Big( C_V \Theta + \frac{U^2}{2} \Big) ' 
+ {\rm \, div \,} \Big( PR\Theta U \frac{x}{r} \Big) 
\\
& 
- \kappa \Delta \Theta 
- (2\mu + \lambda) \Delta \frac{U^2}{2} 
\\
 =
 & \lambda {\rm \, div \,} \Big( \frac{d-1}{r} U^2 \frac{x}{r} \Big) .
\end{split}
\end{equation*}
Let $\varphi \in C_0 ^\infty (\mathbb R^d)$ such that $\varphi (0) = 1$, $\nabla \varphi (0) = 0$, 
${\rm supp } \, \varphi \subset \{ |x| \leq 1  \}$, 
and 
$\varphi_\varepsilon (x) := \varphi (\varepsilon x)$. 
We consider the duality pair of the distribution and the Schwartz function 
for the above equation and $\varphi _{\varepsilon}$. Except for the second term in the 
left-hand side, it follows by the Lebesgue dominated convergence theorem that 
\[
\lim_{\varepsilon \to 0} 
\Big\langle  P \Big( C_V \Theta + \frac{U^2}{2} \Big) , \varphi_\varepsilon  \Big\rangle 
= 
\int_{\mathbb R^d} P \Big( C_V \Theta + \frac{U^2}{2} \Big) dx 
\]
\[
\lim_{\varepsilon \to 0} 
\Big\langle  {\rm \, div \,} \Big( PR\Theta U \frac{x}{r} \Big)  , \varphi_\varepsilon  \Big\rangle 
= \lim_{\varepsilon \to 0} 
\int_{\mathbb R^d} PR\Theta U \frac{x}{r} \cdot \nabla \varphi _\varepsilon dx = 0,
\]
\[
\lim_{\varepsilon \to 0} 
\langle  \Delta \Theta  , \varphi_\varepsilon  \rangle 
= \lim_{\varepsilon \to 0} 
\int_{\mathbb R^d} \Theta \Delta \varphi _\varepsilon dx = 0, 
\qquad \lim_{\varepsilon \to 0} 
\Big\langle  \Delta \frac{U^2}{2}  , \varphi_\varepsilon  \Big\rangle 
= \lim_{\varepsilon \to 0} 
\int_{\mathbb R^d} \nabla \frac{U^2}{2} \cdot \nabla \varphi _\varepsilon dx = 0,
\]
\[
\lim_{\varepsilon \to 0} 
\Big\langle  {\rm \, div \,} \Big( \frac{d-1}{r} U^2 \frac{x}{r} \Big)  , \varphi_\varepsilon  \Big\rangle 
= \lim_{\varepsilon \to 0} 
\int_{\mathbb R^d} \frac{d-1}{r} U^2 \frac{x}{r} \cdot \nabla \varphi _\varepsilon dx = 0. 
\]
For the remained term, it follows by 
\[
 {\rm div \,} 
\Big( P\big( \frac{r}{2}+U\big) \frac{x}{r} \Big) 
=  \frac{d}{2}P
\]
and the integration by parts that 
\[
\begin{split}
& 
\Big\langle  P\Big( \frac{r}{2} + U \Big) \Big( C_V \Theta + \frac{U^2}{2} \Big) ' , 
    \varphi_\varepsilon  
\Big\rangle 
\\
=
& 
- \int_{\mathbb R^d} \frac{d}{2} P\Big( C_V \Theta + \frac{U^2}{2} \Big)  \varphi_\varepsilon ~dx
- \int_{\mathbb R^d} P\Big( \frac{r}{2} + U \Big) \Big( C_V \Theta + \frac{U^2}{2} \Big) 
  \frac{x}{r} \cdot \nabla \varphi _{\varepsilon} ~dx. 
\end{split}
\]
By the Lebesgue dominated convergence theorem, the limit of the first term is 
\[
\lim_{\varepsilon \to 0} (-1) \int_{\mathbb R^d} \frac{d}{2} P\Big( C_V \Theta + \frac{U^2}{2} \Big)  \varphi_\varepsilon ~dx
=  - \int_{\mathbb R^d} \frac{d}{2} P\Big( C_V \Theta + \frac{U^2}{2} \Big) ~dx.
\]
The second term vanishes, since 
\[
\begin{split}
& 
\Big|\int_{\mathbb R^d} P \frac{r}{2} \Big( C_V \Theta + \frac{U^2}{2} \Big) 
  \frac{x}{r} \cdot \nabla \varphi _{\varepsilon} ~dx
\Big| \leq 
\int_{\mathbb R^d} P \frac{ 1 }{2}  \Big( C_V \Theta + \frac{U^2}{2} \Big) 
  \varepsilon r |(\nabla \varphi)(\varepsilon x)| ~dx
  \to 0 ,
\\
& 
\Big| \int_{\mathbb R^d} PU  \Big( C_V \Theta + \frac{U^2}{2} \Big) 
  \frac{x}{r} \cdot \nabla \varphi _{\varepsilon} ~dx
\Big| 
\leq  \varepsilon \|\nabla \varphi \|_{L^\infty} 
    \int_{|x| \leq \varepsilon ^{-1}} PU  \Big( C_V \Theta + \frac{U^2}{2} \Big)  dx
  \to 0 
\end{split}
\]
as $\varepsilon \to 0$. 
Therefore we obtain 
\[
\int_{\mathbb R^d}
 \bigg[P \Big( C_V \Theta + \frac{U^2}{2} \Big) 
 -\frac{d}{2} P\Big( C_V \Theta + \frac{U^2}{2} \Big)  
 \bigg]
 = 0,
\]
which proves $C_V \Theta + \frac{U^2}{2} = 0$ provided that $d \not = 2$.

\section{Large Shrinkers: numerical approach} \label{LSNA}
A natural question in view of Theorem \ref{thm} is whether large shrinkers exist or not. We have investigated this question numerically.
 
\subsection{Description of the simulations} \label{simudescription}
We use the ode45 solver of Matlab to solve \eqref{eq:ODEs_back}. Since the equations are singular, we evaluate the initial data at $r=\delta \ll 1.$ Then, we have to prescribe initial conditions for $P,U,U',\Theta$ and $\Theta'$. Since for $r$ small the backward and forward systems are essentially equivalent, we use the results of the first two author's paper \cite{GI}. In this work, two regimes of initial data consistent with the equations are identified. This provides appropriate choices for our numerical simulations in the present paper. \\
\indent -In the case of cavitating solutions, the choice of initial conditions is straightforward: we take $P(\delta)=P_{\delta}, U(\delta)=-\alpha \delta, U'(0)=-\alpha, \Theta(0)=\Theta_0, \Theta'(0)=0,$ where $P_{\delta}, \alpha$ and $\Theta_0$ are free parameters.
\\
\indent -In the case of smooth solutions, the authors showed in \cite{GI} that once the initial density and temperatures are chosen (say equal to $P_0$ and $\Theta_0)$, the velocity is fixed, and $U(r) = O(r^3).$ To determine an adequate initial condition, we have to find the first term in this asymptotic expansion. It turns out that it depends on further terms in the expansion of $\Theta,$ which leads us to correcting $\Theta$ at the second order. In conclusion, we assume $P(r) = P_0,$ $U(r) = Ar^3$ and $\Theta = \Theta_0 + Br^2$ near the origin, with $A$ and $B$ to be determined. Plugging these ansatz into the integro-differential formulation of the equations (see \cite{GI} for details on the implementation of this procedure), we obtained that as $r \to 0,$
\begin{align*}
U(r) &=  \frac{R P_0 ^2 \Theta_0 C_V}{30 \kappa (2\mu+\lambda+RP_0 \Theta_0)} r^3 + O(r^4), \\
\Theta(r) &= \Theta_0 + \frac{C_V P_0 \Theta_0}{6\kappa} r^2 + O(r^3).
\end{align*}
For the numerical simulations we took $C_V = R = \kappa = \lambda= 1,$ $\mu=2$ and $d=3.$

\subsection{Results}
We now present the results of the numerical simulations described in the previous section.\\
\\
\underline{Cavitating solutions:}
\\
We generically observe that the velocity and temperature diverge to $-\infty,$ and the density blows up to $+\infty$ very fast. This last observation is not surprising given the explicit expression of $P$ in terms of $U,$ obtained by direct integration of the first equation in \eqref{eq:ODEs_back}. Here are a few typical profiles that we obtain (we varied the density from $0.1$ to $1$ and kept the initial velocities and temperatures constant equal to $-0.1$ and $1$ respectively):
\begin{figure}[!h]
\begin{minipage}{.45\textwidth}
\centering
    \includegraphics[width=0.9\linewidth]{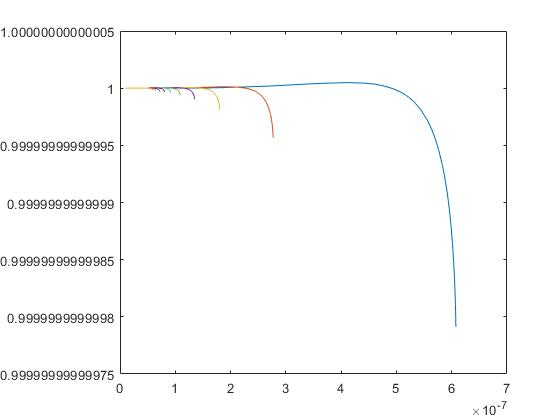}
    \caption{graph of $\Theta$}
    \end{minipage}
\begin{minipage}{.45\textwidth}
\centering
    \includegraphics[width=0.9\linewidth]{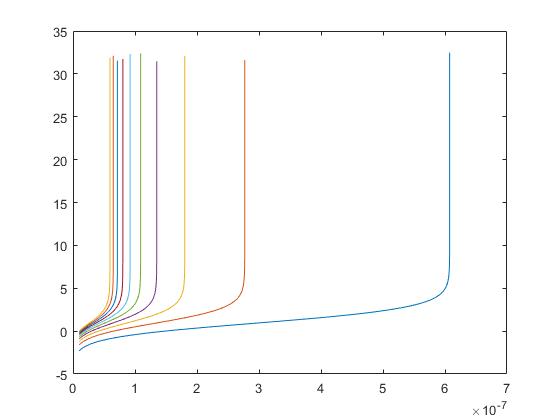}
    \caption{graph of $\log P$}
    \end{minipage}
    \begin{minipage}{.45\textwidth}
\centering
    \includegraphics[width=0.9\linewidth]{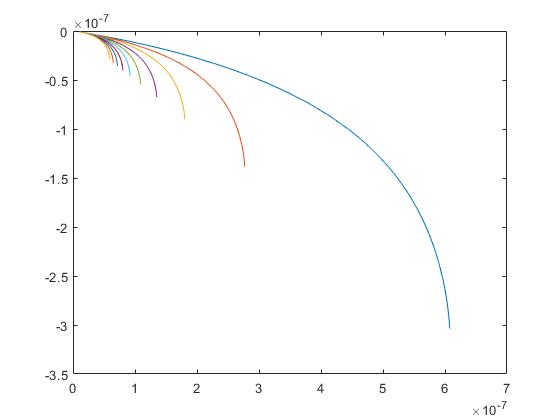}
    \caption{graph $U$}
    \end{minipage}
\end{figure}
\newpage
To test how typical this behavior is, we compute the value of the velocity profile $U(r)$ for some large $r$ and an initial condition $(P_0, \Theta_0, \alpha).$ We represent the corresponding point of coordinates $(P_0, \Theta_0, \alpha)$ in blue if $U(r)$ is negative, and red if it is positive. Indeed if a large shrinker exists, we expect to see a transition in the behavior of $U,$ from a blow-up towards $+\infty$ to a blow-up towards $-\infty.$ The results are the following:
\begin{figure}[!h]
\centering
\includegraphics[width=0.55\linewidth]{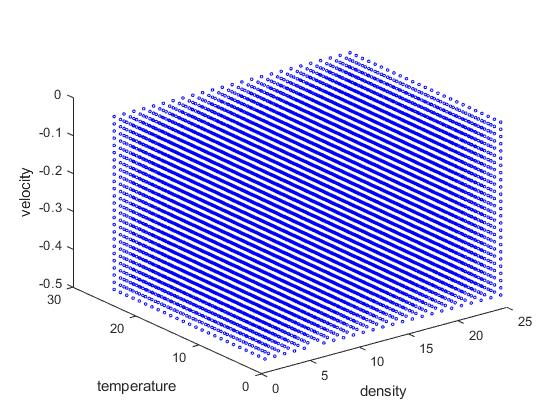}
\caption{$\delta=10^{-3}$}
\end{figure}
\begin{figure}[!h]
\centering
\includegraphics[width=0.55\linewidth]{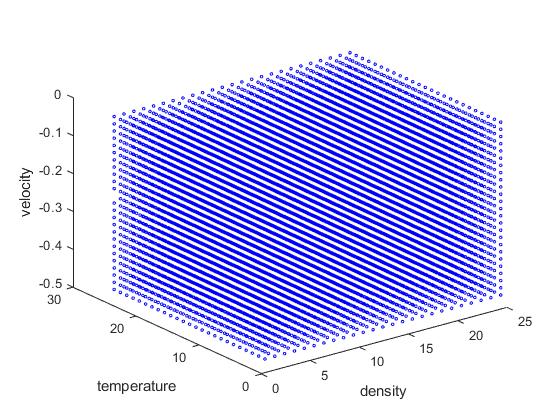}
\caption{$\delta=10^{-6}$}
\end{figure}
\newpage
\noindent Note that in both cases, there are no red points.
\\
\\
\underline{Smooth solutions:}\\
In this case, we observe the following typical behavior: $U$ and $\Theta$ blow-up to infinity, and $P$ converges to 0. The derivative of $U$ also blows-up to infinity, at a seemingly faster rate than $U.$ Below are a few typical graphs of solutions that we obtained (they all correspond to $P_0$ ranging from 0.3 to 3, and initial temperature equal to 1). 
\begin{figure}[!h]
\begin{minipage}{.45\textwidth}
\centering
    \includegraphics[width=0.9\linewidth]{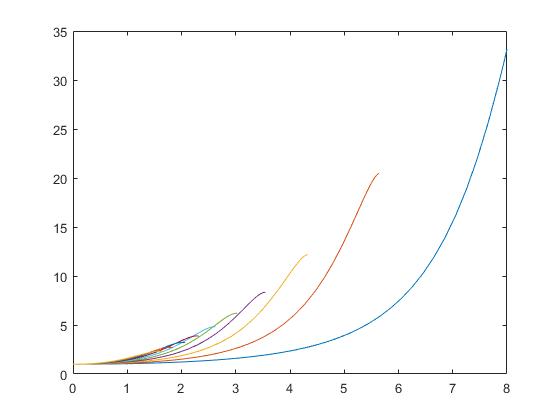}
    \caption{Temperature}
    \end{minipage}
\begin{minipage}{.45\textwidth}
\centering
    \includegraphics[width=0.9\linewidth]{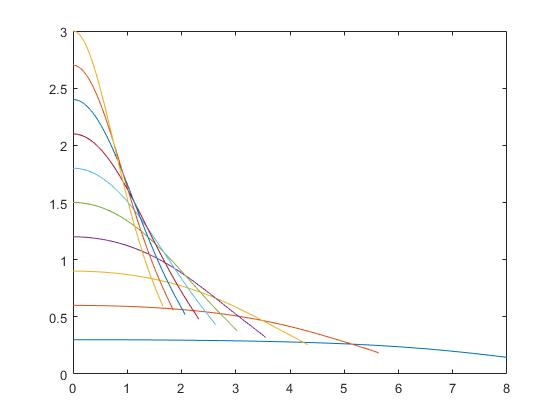}
    \caption{Density}
    \end{minipage}
\end{figure}
\begin{figure}[!h]
    \begin{minipage}{.45\textwidth}
\centering
    \includegraphics[width=0.9\linewidth]{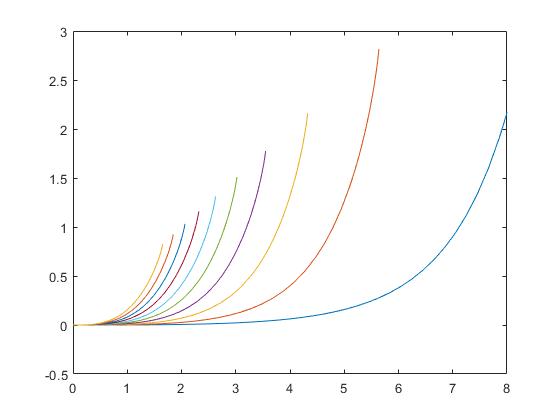}
    \caption{Velocity}
    \end{minipage}
    \begin{minipage}{.45\textwidth}
    \centering
    \includegraphics[width=0.9\linewidth]{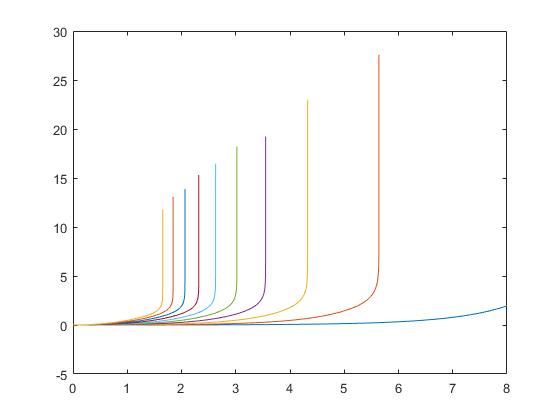}
    \caption{Derivative of velocity}
    \end{minipage}
\end{figure}
\newpage
\noindent As we did for cavitating solutions, we test how universal this behavior is: we compute values of $U$ for $r$ large and $P_0,$ $\Theta_0$ between 0 and 40. If the value obtained is positive, we represent the point of coordinates $(P_0,\Theta_0)$ in blue on the picture below. If it is negative, we represent it in red. In this case the program occasionally returns an error. In this case the corresponding point is green. We obtain the following (for various values of $\delta$):
\begin{figure}[!h]
\centering
\includegraphics[width=0.75\linewidth]{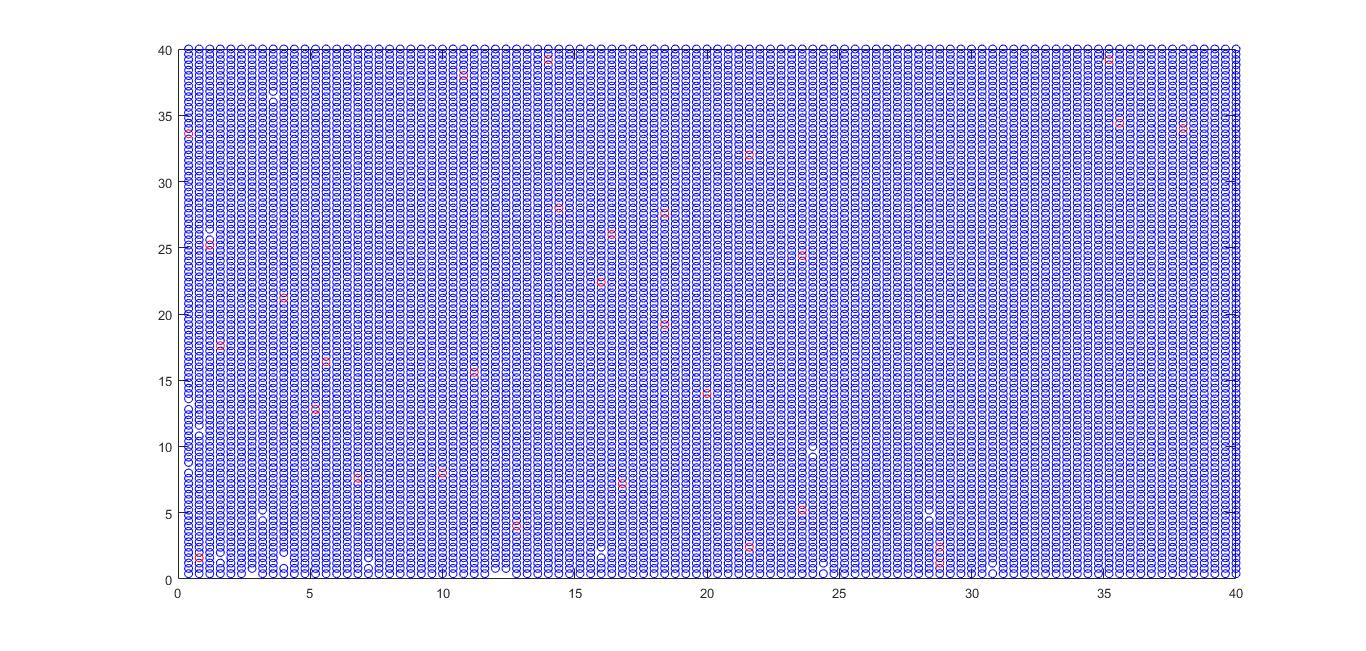}
\caption{$\delta=10^{-5}$}
\end{figure}
\begin{figure}[!h]
\centering
\includegraphics[width=0.75\linewidth]{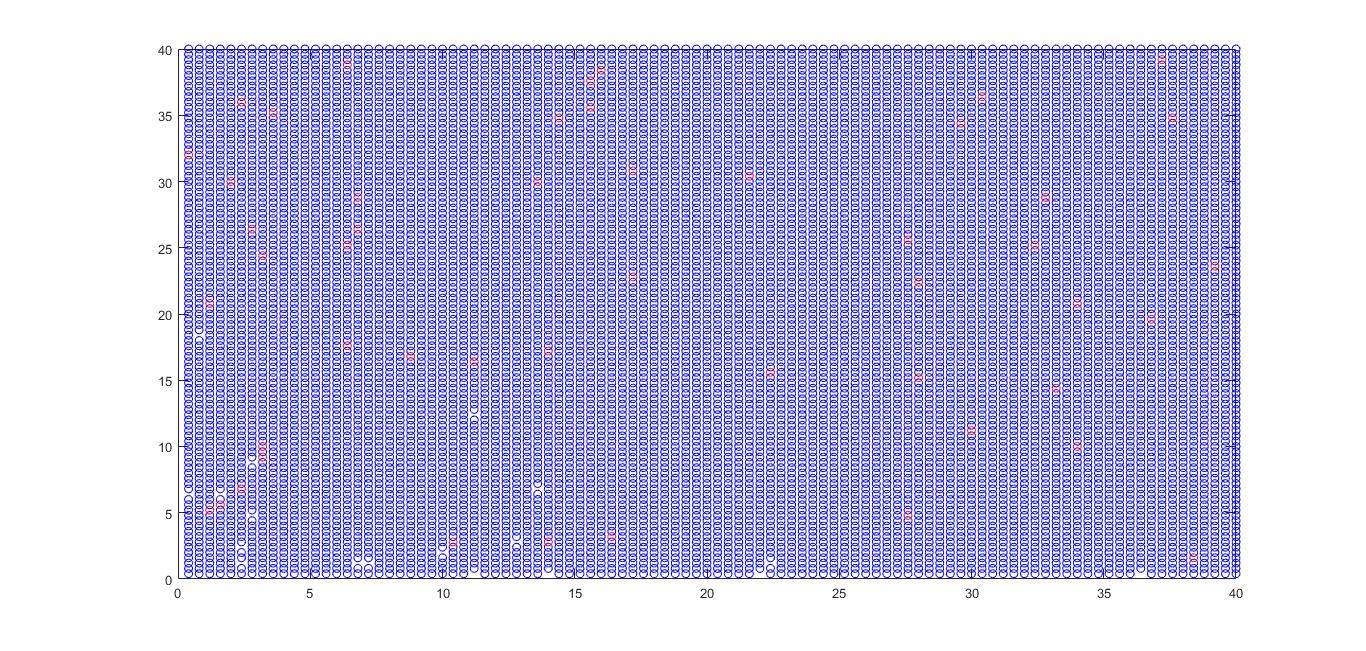}
\caption{$\delta=10^{-8}$}
\end{figure}
\begin{figure}[!h]
\includegraphics[width=0.75\linewidth]{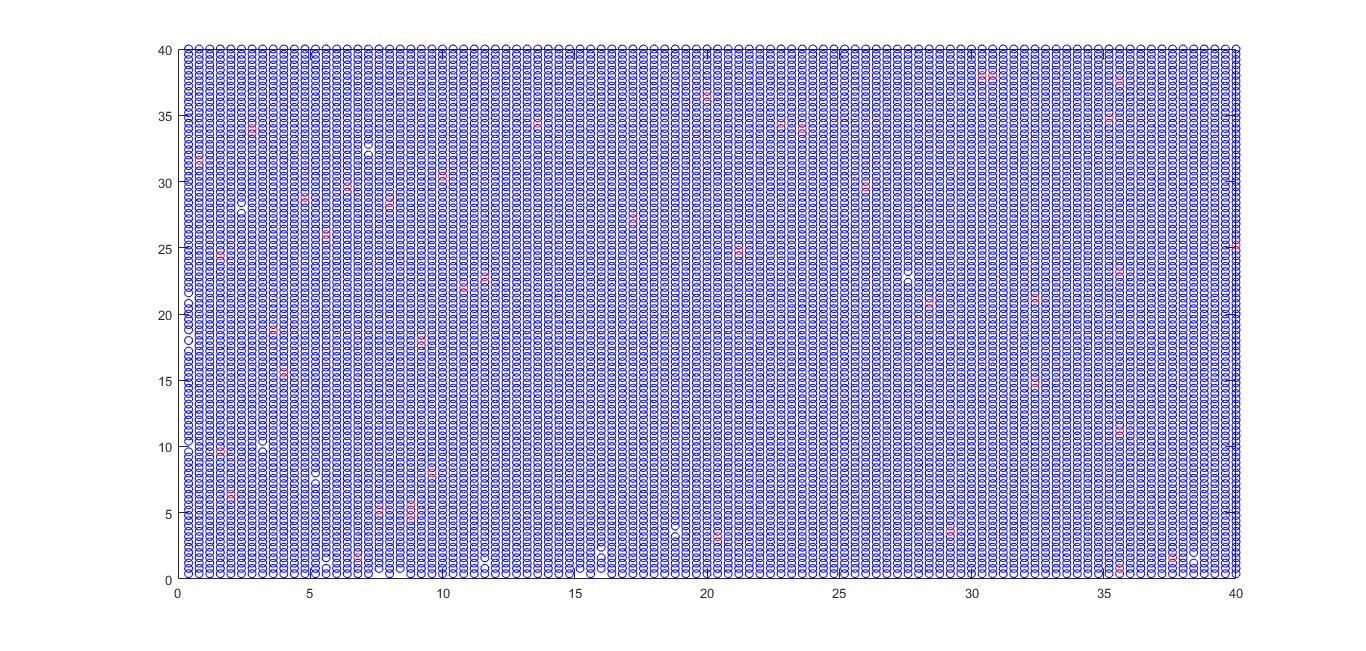}
\caption{$\delta=10^{-9}$}
\end{figure}
\newpage
\noindent \underline{Conclusion:} \\
In both cases the blow up behavior identified appears generic. Note that results in the smooth case are slightly less decisive since there are a few anomalous red points. However they seem unstable numerically (their position changes as $\delta$ varies). Their presence may be due to the fact that the initial data for the velocity is obtained implicitly in this case, and that in fact we only consider the first term in the expansion of the profile (see subsection \ref{simudescription}). These simulations support the claim that there are no large radially symmetric shrinkers for \eqref{eq:cNS}. 


\end{document}